\documentclass{amsart}
\usepackage{tikz}
\usetikzlibrary{calc,through,backgrounds}
\usepackage{xcolor}
\usepackage{multicol}
\usepackage{bbm} 
\usepackage{bm,amssymb,amsthm,amsfonts,amsmath,latexsym,cite,color,psfrag,graphicx,ifpdf,pgfkeys,pgfopts,xcolor,extarrows}
\usepackage{todonotes}

\newtheorem{theo}{Theorem}[section]
\newtheorem{defi}[theo]{Definition}
\newtheorem{exam}[theo]{Example}
\newtheorem{lem} [theo]{Lemma}
\newtheorem{coro}[theo]{Corollary}
\newtheorem{prop}[theo]{Proposition}

\newtheorem{re}[theo]{Remark}
\newtheorem{conj}[theo]{Conjecture}

\allowdisplaybreaks
\def\P{   \mathcal{P}  }
\def\rk{   {\texttt{rk}}  }
\def\SM{   { \texttt{SM} } }

\makeatletter \@addtoreset{equation}{section}
\@addtoreset{theo}{}\makeatother

\usepackage{hyperref}

\DeclareMathAlphabet{\amathbb}{U}{bbold}{m}{n}

\newcommand{\indicator}[1]{\overline{\amathbb{1}}\!\left({#1}\right)}

\def\qed{\hfill \rule{4pt}{7pt}}
\def\pf{\noindent {\it Proof.} }

\begin{document}
\begin{center}
{\large\bf  VALUATIVE INVARIANTS OF CATALAN MATROIDS }
			
\vskip 4mm
{\small YIMING CHEN,  YAO LI, AND MING YAO}
	
\end{center}
	
\noindent{\textsc{Abstract.}}
We decompose the indicator function of each $(a, b)$-Catalan matroid polytope as a weighted sum of indicator function of matroid polytopes that correspond to direct sums of uniform matroids. Catalan matroids lie in the interior of the convex hull of direct sums of uniform matroids in the polytope of all matroids introduced by Ferroni and Fink. Moreover, we describe combinatorially the coefficients of the convex combination of direct sums of uniform matroid corresponding to an $(a,b)$-Catalan matroid. In particular, this allows us to derive explicit formulas for arbitrary valuative invariants of $(a, b)$-Catalan matroids. Among other applications, we prove that $(a,b)$-Catalan matroids are Ehrhart positive, and we find formulas for the Kazhdan--Lusztig invariants of these matroids.

\section{Introduction}

Catalan matroids, a subfamily of Schubert matroids, were first studied by Ardila \cite{Ard-1} and independently by Bonin, de Mier, and Noy \cite{Bon}. The latter authors also introduced $(a,b)$-Catalan matroids.
For $a,b\ge1$, let $U_n^{a,b}:=U_{nb,n(a+b)}$ denote the uniform matroid on $[n(a+b)]$ with bases $\{B\colon B\subseteq [n(a+b)],|B|=nb\}$, and let $C_n^{a,b}$ be the   $(a,b)$-Catalan matroid on $[n(a+b)]$.  When $a=b=1$, write $U_n:=U_{n}^{1,1}$ and $C_n:=C_n^{1,1}$ for short.  

Suppose that $M=([n], \mathcal{B})$ is a matroid with ground set $[n]:=\{1,2,\ldots,n\}$  and base set $\mathcal{B}$. The  matroid polytope $\P(M)$ of $M$ is defined as the convex hull
$
\P(M)=\mathrm{conv}\{e_B\colon B\in \mathcal{B}\}, $
where $e_B=e_{b_1}+\cdots+e_{b_k}$ for $B=\{b_1,\ldots,b_k\}\subseteq[n]$, and $e_i$ is  the $i$-th standard basis of $\mathbb{R}^n$. 
 Let $X$ be a subset of $\mathbb{R}^n$, its indicator function $\indicator{X}: \mathbb{R}^n \rightarrow \mathbb{Z}$ is defined by
$$
\indicator{X}(x)= \begin{cases}1 & x \in X, \\ 0 & \text {otherwise.}\end{cases}
$$
The polytope of all matroids\cite{F-allmatroids} arises from Derksen and Fink's construction\cite{Schubert-expansion}, which shows that the indicator function of any matroid polytope can be expressed as an integer linear combination of Schubert matroid polytope indicators. In other words, each matroid $M$ on $[n]$ of rank $k$ admits a unique Schubert expansion
\begin{align}\label{pm}
\indicator{\P(M)} \;=\;\sum_{\text{S}} c_{\text{S}}\,\indicator{\P(\text{S})},
\end{align} 
where the sum runs over Schubert matroids on $[n]$ of rank $k$.  Ferroni and Fink \cite{F-allmatroids} characterize the coefficients $c_{\text{S}}$  in \eqref{pm}, which  form a vector $p_{M}$ and provides a natural coordinate for the matroid $M$.

 The convex hull of the coordinate vectors 
  $p_{M}$ for all matroids $M$ on $[n]$ of rank $k$ defines the labeled polytope $\overline{\Omega}_{k,n}$. By identifying isomorphic matroids and taking the natural projection, one obtains the unlabeled polytope $\Omega_{k,n}$. We denote the isomorphism class of a matroid $M$ by $[M]$, and the coordinate vector of $[M]$ by $p_{[M]}$. These polytopes serve as parameter spaces for valuative invariants and provide a geometric framework for investigating non-negativity conjectures in matroid theory.
   
Ferroni and Fink \cite{F-allmatroids} showed that in $\overline{\Omega}_{k,n}$, every matroid corresponds to a vertex. For $\Omega_{k,n}$, however, the situation is more delicate: distinct matroid isomorphism classes may collapse to the same point, and not every isomorphism class corresponds to a vertex.

\begin{defi}[\cite{F-allmatroids}, Definition 4.2]
    A matroid $M$ of rank $k$ on $n$ elements is said to be \emph{extremal} if its coordinate vector $p_{[M]}$ is a vertex of the polytope $\Omega_{k,n}$.
\end{defi}
 
In this paper, we show that  
\begin{theo}\label{Theo_nonextrme}
    Catalan matroids are non-extremal matroids. Moreover, the coordinate vector $p_{[C_{n}]} \in \Omega_{n, 2 n}$ of the Catalan matroid $C_{n}$ lies in the interior of the convex hull of the coordinate vectors of direct sum of uniform matroids.
\end{theo}

To prove Theorem \ref{Theo_nonextrme}, we first realize the uniform matroid polytope as a union of $(a,b)$-Catalan matroid polytopes via a matroid subdivision, then decompose the matroid polytope of $C_n^{a,b}$ into a direct product of (smaller) uniform matroids polytopes. This enables us to derive the closed-form formula \eqref{main-formula} for  valuative invariants of $C_n^{a,b}$, also yields nontrivial applications of these invariants.

Let $A$ be an abelian group and $\overline{\mathcal{M}}_{k,n}:=\{\text{matroids on $[n]$ of rank $k$}\}$. A map $f:\overline{\mathcal{M}}_{k,n}\to A$ is a \emph{valuation} if $\sum_{i=1}^{k}a_{i}f(M_{i})=0$ whenever $\sum_{i=1}^{k}a_{i}\indicator{\mathcal{P}(M_{i})}$ is the identically zero function for all coefficients $a_{i}\in \mathbb{Z}$.
Let $\mathcal{M}_{k,n}$ be the set of isomorphism class of   matroids on $[n]$ of rank $k$. A {\it valuative invariant} is a function $f:\mathcal{M}_{k,n}\to A$ such that its induced function   on $\overline{\mathcal{M}}_{k,n}$ is a valuation.

A \emph{matroid subdivision} of $\mathcal{P}(M)$ is a polyhedral complex whose maximal cells are all matroid  polytopes and whose union is $\mathcal{P}(M)$.
Given   a subdivision of $\mathcal{P}(M)$ with maximal cells $\mathcal{P}(M_1), \dots, \mathcal{P}(M_s)$, any valuative function $f$ satisfies the inclusion–exclusion principle
\[
f(M) = \sum_{J \subseteq [s]} (-1)^{|J|-1} f\left(\bigcap_{j \in J} M_j\right),
\]
where $\bigcap_{j \in J} M_j$ is the matroid whose base polytope is  $\bigcap_{j \in J} \mathcal{P}(M_j)$. When the intersection is empty, we adopt the convention that $f(\bigcap_{j \in J} M_j) := 0$.

Recall that a sequence of positive integers $\lambda=
(\lambda_1,\ldots,\lambda_k)$ is called a partition of 
$n$, denoted as $\lambda\vdash n$, if $\lambda_1\ge\cdots
\ge\lambda_k\ge1$ and $\lambda_1+\cdots+\lambda_k=n$. Let
\[z_{\lambda}=\prod_{i=1}^ni^{\alpha_i}\alpha_i!,\]
where $\alpha_i$ is  the number of appearances of $i$ in $\lambda$.
A sequence of positive integers $\gamma=(\gamma_1,\ldots,\gamma_k)$ is called a composition of $n$,  if $\gamma_1+\cdots+\gamma_k=n$.  
For any valuative function, we define
\[ f(U_{\gamma}^{a,b}):=\prod_{i=1}^k f(U_{\gamma_i}^{a,b}).\]

\begin{coro}
\label{CoroVI}
Let $f$ be a valuative invariant, we have
    \begin{equation}\label{main-formula}
    f(C_n^{a,b})=\sum_{\lambda\vdash n} \frac{1}{z_{\lambda}} f(U_{\lambda}^{a,b})
    \end{equation}
\end{coro}

As applications of Corollary \ref{CoroVI}, we derive the Ehrhart polynomials of $(a,b)$-Catalan matroids and obtain their Ehrhart positivity. Furthermore, we compute the volume, Tutte polynomial, (inverse) Kazhdan–Lusztig polynomial, Z-polynomial and Whitney polynomial of $(a,b)$-Catalan matroids.

\section{Preliminaries}

A {\it matroid} is a pair $M=(E, \mathcal{I})$ consisting of a finite set $E$,  and a collection $\mathcal{I}$  of subsets of $E$, satisfies:
(i) $\emptyset \in \mathcal{I}$;
(ii) If $J\in \mathcal{I}$ and $I\subseteq J$, then $I\in \mathcal{I}$;
and (iii)  If $I, J\in \mathcal{I}$ and $|I|<|J|$, then there exists $j\in J\setminus I$
such that $I\cup \{j\}\in \mathcal{I}$. 
By (ii), a matroid $M$ is determined by the collection $\mathcal{B}$  of maximal independent sets, called the bases of $M$. By (iii),  all the bases have the same cardinality, called the rank of $M$, denoted as ${\rk}(M)$. So we can  write $M=(E,\mathcal{B} )$. The {\it dual} of $M$ is a 
matroid $M^*=(E,\mathcal{B}^*)$ where $\mathcal{B}^*=\{E\setminus B\colon B\in \mathcal{B}(M)\}.$
The {\it rank  function} ${\rk}_M: 2^E\rightarrow \mathbb{Z}$  of $M$  is  defined as
 \[{\rk}_M(T)=\max\{|T\cap B|\colon B\in \mathcal{B}\},\ \ \ \text{for $T\subseteq E$}.\]
The {\it closure} of a set $A\subseteq E$ is  the biggest set that contains $A$ and has the same rank. A set is a {\it  flat} if its closure is same
as itself. A set $A\subseteq E$ is called {\it separable} if 
there are disjoint subsets $R$ and $T$ such that $A=R\cup T$ and ${\rk}_M(A) = {\rk}_M(R)+{\rk}_M(T).$ Otherwise, $A$ is called {\it inseparable}.

Let  $[n]:=\{1,2,\ldots,n\}$ and $S$ be a subset of $[n]$. The {\it Schubert matroid} ${\SM}_n(S)$ is the matroid with ground set $[n]$ and  bases
\[\{T\subseteq [n]\colon T\leq S\},\]
where $T\leq S$ means that: $|T|=|S|$ and
the $i$-th smallest element of $T$ is less than or equal to that of $S$ for $1\leq i\leq |T|$.
The {\it indicator vector}  of $S$ is the 0-1 vector $\mathbb{I}(S)=(i_1,\ldots,i_n)$, where $i_j=1$ if $j\in S$, and 0 otherwise. 
Without loss of generality, we can assume that $1\not\in S$ and the largest element of $S$ is $n$. Thus $i_1=0$ and $i_n=1$.  
For simplicity, write $\mathbb{I}(S)=(0^{r_1},1^{r_2},\ldots,0^{r_{2m-1}},1^{r_{2m}})$,
where $0^{r_1}$ represents $r_1$ copies of 0's, $1^{r_2}$ represents $r_2$ copies of 1's, etc. Then $S$ can be written   as a  sequence of positive integers $r(S)=(r_1,r_2,\ldots,r_{2m})$ of length $2m$.
It is easy to see that given such an integer sequence $r$, there is a unique set $S$  whose indicator vector $\mathbb{I}(S)$ can be written   in this way. 
For example, let $S=\{2,5,7,8,9\}\subseteq[9]$, then $\mathbb{I}(S)=(0,1,0^2,1,0,1^3)$ and $r(S)=(1,1,2,1,1,3)$.

Let $r(S)=(n-k,k),(1,k-1,n-k-1,1),(k-1,1,1,n-k-1)$, we obtain the uniform matroids $U_{k,n}$ \cite{Fer}, the minimal matroids $T_{k,n}$ \cite{Fer2}, and the sparse paving Schubert matroids ${\texttt{Sp}}_{k,n}$ \cite{FL}, respectively. For any positive integers $a,b,c,d$, $r(S)=(a,b,c,d)$ corresponds to the notched rectangle matroids  \cite{FL}, which are conjectured to be Ehrhart positive. In particular,   $r(S)=(a,b,c,1)$ corresponds to the {\it panhandle matroids} ${\texttt{Pan}}_{b+1,a+b,n}$ \cite{panhandle,panhandleposi}, where $n=a+b+c+1$.  Moreover, let $r(S)=(\overbrace{a,b,a,b,\ldots,a,b}^{2n})$,  we obtain the {\it $(a,b)$-Catalan matroids}  $C_n^{a,b}$ introduced by Bonin,  de Mier and  Noy \cite[Definition 3.7]{Bon}. Notice that $C_n^{a,b}$ is a matroid on 
$[n(a+b)]$. When $a=b=1$,  $C_n^{1,1}$ is the $n$-th Catalan matroid $C_n$ as first studied by Ardila \cite{Ard-1}.

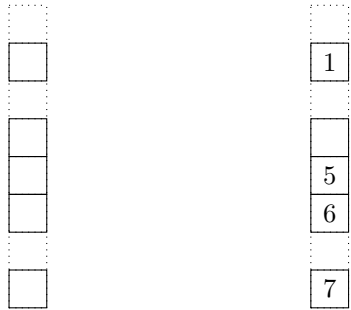
\begin{figure}[ht]
\begin{center}
\begin{tikzpicture}

\def\rectanglepath{-- +(5mm,0mm) -- +(5mm,5mm) -- +(0mm,5mm) -- cycle}

\draw [step=5mm,dotted] (0mm,0mm) grid (5mm,40mm);
\draw (0mm,0mm) \rectanglepath;
\draw (0mm,10mm) \rectanglepath;
\draw (0mm,15mm) \rectanglepath;
\draw (0mm,20mm) \rectanglepath;
\draw (0mm,30mm) \rectanglepath;

\draw [step=5mm,dotted,shift={(40mm,0mm)}] (0mm,0mm) grid (5mm,40mm);

\draw[shift={(40mm,0mm)}] (0mm,0mm) \rectanglepath;
\draw[shift={(40mm,0mm)}] (0mm,10mm) \rectanglepath;
\draw[shift={(40mm,0mm)}] (0mm,15mm)\rectanglepath;
\draw[shift={(40mm,0mm)}] (0mm,20mm) \rectanglepath;
\draw[shift={(40mm,0mm)}] (0mm,30mm) \rectanglepath;
\node [shift={(40mm,0mm)}]  at (2.5mm,2.5mm) {$7$};
\node [shift={(40mm,0mm)}]  at (2.5mm,12.5mm) {$6$};
\node [shift={(40mm,0mm)}]  at (2.5mm,32.5mm) {$1$};
\node [shift={(40mm,0mm)}]  at (2.5mm,17.5mm) {$5$};

\end{tikzpicture}
\end{center}
\vspace{-4mm}
\caption{The filling $\mathcal{F}_{S}(w)$ for $w=67185$.}
\label{diagram}
\end{figure}

It is well known that the matroid polytope $\P(M)$ defined in \eqref{pm}   is a generalized permutohedron parametrized by the rank function of $M$, see, for example, Fink, M\'esz\'aros and St.$\,$Dizier \cite{Fin}.  More precisely,
\begin{equation}\label{DE}
\P(M)=\left\{x\in \mathbb{R}_{\ge0}^n\colon \sum_{i\in [n]}  x_i={\rk}_M([n])
\ \  \text{and}\ \ \sum_{i\in T}x_i\leq  {\rk}_M(T)\ \ \text{for $T\subsetneq [n]$}\right\}.
\end{equation}
Feichtner and  Sturmfels \cite{FS} provided a minimal system of inequalities for all matroid polytopes \eqref{DE}.
In order to give a minimal system of inequalities for Schubert matroid  polytopes, we recall the  
combinatorial algorithm of computing the rank function $\rk_S(T)$ of a Schubert matroid $\SM_n(S)$ provided by Fan and Guo \cite[Theorem 3.3]{FL1}.
Given a subset $S=\{s_1,\ldots,s_k\}$ on $[n]$, view $S$  as   $n\times 1$ grids such that there is a box in row $j$ for $j\in S$, and the other rows are empty. The left diagram in Figure \ref{diagram} is an illustration for $S=\{2,4,5,6,8\} \subseteq[8]$.

For a permutation
$w=w_1w_2\cdots w_t$ of elements of   $T\subseteq[n]$,
define $\mathcal{F}_{S}(w)$ to be the filling
of $S$ with $w$ as follows.
 For $j$ from $1$  to $t$, put $w_j$ into the
 first (from top to bottom) empty box whose row index is larger than
 or equal to $w_j$. If  there are no such empty boxes, then $w_j$ does
 not appear in the filling and skip to $w_{j+1}$.
 Let $|\mathcal{F}_{S}(w)|$ denote the number of non-empty boxes of $S$ after filling the elements of $w$. The right diagram in Figure \ref{diagram} is the filling $\mathcal{F}_{S}(w)$ of $S=\{2,4,5,6,8\}$ with $w=67185$, where $|\mathcal{F}_{S}(w)|=4$.

\begin{theo}[\cite{FL1}]\label{Main-1}
 Let $w=w_1w_2\cdots w_t$ be any   permutation of elements of $T\subseteq[n]$.
Then
\begin{equation}\label{ff}
\rk_S(T)=|\mathcal{F}_{S}(w)|.
\end{equation}
\end{theo}

We have the following minimal system of inequalities for Schubert matroid polytopes.

\begin{lem}\label{plty}
For any $S=\{s_1,s_2,\ldots,s_k\}\subseteq[n]$, we have
\begin{align}\label{C}
    \mathcal{P}({\SM}_n(S))=\left\{x\in \mathbb{R}^{n}_{\ge0}:\sum_{i\in[n]}x_i=k, \ \sum_{i=1}^{s_j}x_i\ge j ,\ \text{for}\ j=1,2,\ldots,k\right\}.
\end{align}
\end{lem}
This description follows from the fact that the inequalities defining a matroid polytope are determined by its cyclic flats. 
Moreover, the cyclic flats of a Schubert matroid form a chain, which yields the inequalities in~\eqref{C}.

\section{The Interior Position of Catalan Matroids}

Recall that the Catalan matroid $C_n$ corresponds to 
$r(S)=(\overbrace{1,1,\ldots,1,1}^{2n})$.
By Lemma \ref{plty}, we have
\begin{align}\label{CC}
    \mathcal{P}(C_n)=\left\{x\in \mathbb{R}^{2n}_{\ge0}:\sum_{i\in[2n]}x_i=n, \sum_{i\in[2j]}x_i\ge j ,\ \text{for}\ j=1,2,\ldots,n-1\right\}.
\end{align}
For $m=0,1,\ldots,n-1$,   define
\begin{align}\label{pcnm}
\mathcal{P}(C_{n,m}):=\left\{x\in \mathbb{R}^{2n}_{\ge0}:\sum_{i\in[2n]}x_i=n,  \sum_{i\in[2j]+2m}x_i \ge j,\ \text{for}\  j=1,2,\ldots,n-1\right\},
\end{align}
where, for simplicity, $x_i$ is understood as $x_{i-2n}$ if $i>2n$.
Clearly, $\mathcal{P}(C_{n})=\mathcal{P}(C_{n,0})$.
Moreover, it is well known that
\[
\mathcal{P}(U_n)=\left\{x\in \mathbb{R}^{2n}_{\ge0}:\sum_{i\in[2n]}x_i=n\right\}.
\]

\begin{prop} We have
\begin{align}\label{u2c}
\mathcal{P}(U_n)=\bigcup_{j=1}^{n} \mathcal{P}(C_{n,j-1}).
\end{align}
\end{prop}

\pf
It is obvious that $\mathcal{P}(C_{n,j-1})\subseteq\mathcal{P}(U_n)$ for any $j=1,\ldots,n$. Let $(x_1,\ldots,x_{2n})\in \mathcal{P}(U_n),$ i.e., $x_1+\cdots+x_{2n}=n$,
we aim to show that 
$(x_1,\ldots,x_{2n})\in\bigcup_{j=1}^{n} \mathcal{P}(C_{n,j-1}).
$

For $i=1,2,\ldots,n$, let $c_i=x_{2i-1}+x_{2i}-1$. If $c_1+\cdots+c_k\ge0$  for all $k=1,2,\ldots,n$, then it is clear that $(x_1,\ldots,x_{2n})\in\mathcal{P}(C_{n})=\mathcal{P}(C_{n,0})$. Suppose that there exists some index $k$ such that $c_1+\cdots+c_k<0$. Let $m=\min\{c_1,c_1+c_2,c_1+c_2+c_3,\ldots,c_1+\cdots+c_n\}$. Then $m<0.$ Choose the smallest index $l$ such that $c_1+\cdots+c_l=m$. We claim that $(x_1,\ldots,x_{2n})\in\mathcal{P}(C_{n,l})$.  

Actually,  the sequence $(c_1,c_2,\ldots,c_n)$ can be viewed as a  path from $(0,0)$ to $(n,0)$ such that each $c_i$ is a step from $(i-1,c_{i-1})$ to $(i,c_i)$, where $c_0=0.$  By cutting off the first $l$ steps of $(c_1,\ldots,c_n)$ and appending them to the end, we obtain $(c_{l+1},c_{l+2},\ldots,c_n,c_1,\ldots,c_l)$. It is easy to see that $(c_{l+1},c_{l+2},\ldots,c_n,c_1,\ldots,c_l)$ is a  path whose steps never go below the $x$-axis. That is, $(c_{l+1},c_{l+2},\ldots,c_n,c_1,\ldots,c_l)\in \mathcal{P}(C_{n,0}),$ which implies  $(x_1,\ldots,x_{2n})\in\mathcal{P}(C_{n,l})$.
\qed

For any subset $A\subseteq[n]$, denote
\begin{align*}
Q(A)=\bigcap_{j\in A} \mathcal{P}(C_{n,j-1}).
\end{align*}
By \eqref{u2c} and the inclusion-exclusion principle, we have
\begin{align}\label{incexc}
	\indicator{U_n}=\sum_{A\subseteq[n]}(-1)^{|A|-1}\cdot \indicator{Q(A)}.
\end{align}
Assume that
  $A=\{a_1,\ldots,a_r\}\subseteq[n]$ arranged increasingly, let $\lambda(A)$ be the partition obtained by arranging the  integers 
\[a_1,a_2-a_1,\ldots,a_r-a_{r-1},n-a_r+a_1\]
decreasingly. Actually, if we put  $1,2,\ldots,n$ on a circle, i.e., $1$  and $n$ are viewed as adjacent  to each other, then the elements of $\lambda(A)$  are gaps between the $a_i$'s.  For example, let $A=\{1,2,4,5\}\subseteq[7]$, then $\lambda(A)=(3,2,1,1)$.

\begin{lem}\label{lm1}
Let $\lambda\vdash n$ be a partition of $n$ and $A\subseteq[n]$   such that $\lambda(A)=\lambda$. Then 
\[
\indicator{Q(A)}=\indicator{C_\lambda}.
\]
\end{lem}

\pf Suppose that $A=\{a_1,\ldots,a_r\}_\le$ and by convention let $a_0=0$.  
	By \eqref{pcnm}, we have
	\begin{align}
		Q(A)&=\bigcap_{i=1}^r
		\left\{x\in \mathbb{R}^{2n}_{\ge0}:\sum_{k\in[2n]}x_k=n,  \sum_{k\in[2j]+2a_i-2}x_k\ge j,\ \text{for}\  1\le j\le n-1\right\}\nonumber\\
		&=\left\{x\in \mathbb{R}^{2n}_{\ge0}:\sum_{k\in[2n]}x_k=n,  \sum_{k\in[2j]+2a_i-2}x_k\ge j,\ \text{for}\ 1\le i\le r \ \text{and}\ 1\le j\le n-1 \right\}. \label{Ine}
	\end{align}
Notice that there are $r(n-1)$ inequalities in \eqref{Ine}. 
For any $1\le i\le r-1$, since
\[
\sum_{k\in[2j]+2a_i-2}x_k\ge j,\ \text{for}\  1\le j\le n-1, 
\]
let $j=a_{i+1}-a_{i}$,
we have 
\begin{align}\label{bds1}
x_{2a_i-1}+x_{2a_i}+\cdots+x_{2a_{i+1}-2}\ge a_{i+1}-a_i.
\end{align}
By letting $j=n-a_{i+1}+a_i$ in
\[
\sum_{k\in[2j]+2a_{i+1}-2}x_k\ge j,\ \text{for}\  1\le j\le n-1,
\]
we get
\[
x_{2a_{i+1}-1}+\cdots+x_{2n+2a_i-2}=x_{2a_{i+1}-1}+\cdots+x_{2n}+x_1+\cdots+x_{2a_i-2}\ge n-a_{i+1}+a_i,
\] 
which is equivalent to
\begin{align}\label{bds2}
x_{2a_i-1}+\cdots+x_{2a_{i+1}-2}\le a_{i+1}-a_i.
\end{align}
Combining \eqref{bds1} and \eqref{bds2}, we obtain that 
\begin{align}\label{fst1}
x_{2a_i-1}+\cdots+x_{2a_{i+1}-2}= a_{i+1}-a_i.
\end{align}

Moreover, by putting $j=n-a_r+a_1$ in
\[
\sum_{k\in[2j]+2a_{r}-2}x_k\ge j,\ \text{for}\  1\le j\le n-1,
\]
we have 
\begin{align}\label{bds3}
x_{2a_{r}-1}+\cdots+x_{2n+2a_1-2}=x_{2a_{r}-1}+\cdots+x_{2n}+x_1+\cdots+x_{2a_1-2}\ge n-a_{r}+a_1.
\end{align}
By plugging $j=a_r-a_1$ into 
\[
\sum_{k\in[2j]+2a_{1}-2}x_k\ge j,\ \text{for}\  1\le j\le n-1,
\]
we find
\[
x_{2a_{1}-1}+\cdots+x_{2a_{r}-2}\ge a_{r}-a_1,
\]
which is equivalent to
\begin{align}\label{bds4}
x_1+\cdots+x_{2a_1-2}+x_{2a_{r}-1}+\cdots+x_{2n}\le n-a_{r}+a_1.
\end{align}
Combining \eqref{bds3} and \eqref{bds4}, we obtain that 
\begin{align}
x_1+\cdots+x_{2a_1-2}+x_{2a_{r}-1}+\cdots+x_{2n}= n-a_{r}+a_1.
\end{align}
Therefore, we have
\begin{align}\label{spl}
\left\{
 \begin{array}{lll}
x_{2a_i-1}+\cdots+x_{2a_{i+1}-2}=a_{i+1}-a_{i},\quad  \hbox{for $1\le i\le r-1$,} \\
x_1+\cdots+x_{2a_1-2}+x_{2a_{r}-1}+\cdots+x_{2n}= n-a_{r}+a_1, \\
x_{2a_i-1}+\cdots+x_{2a_{i}-2+2j}\ge j,\quad \hbox{for $1\le j\le n-1$ \ and \ $1\le i\le r$.}
\end{array}
\right.
\end{align}
On the other hand, it's easy to see that \eqref{spl} implies  \eqref{Ine}, so we conclude that \eqref{Ine} is equivalent to \eqref{spl}, namely, $Q(A)$ is the set of nonnegative solutions of  equation \eqref{spl}.

We proceed to show that the set of nonnegative solutions of equation \eqref{spl} is isomorphic to the direct sum of some Catalan matroid polytopes.
For $1\le i\le r-1$, by \eqref{C}, 
\[
\mathcal{P}(C_{a_{i+1}-a_i})=\left\{x\in \mathbb{R}^{2n}_{\ge0}: \sum_{k\in[2(a_{i+1}-a_i)]}x_k=a_{i+1}-a_{i},  \sum_{k\in[2j]}x_k\ge j ,\ \text{for}\ 1\le j\le a_{i+1}-a_i-1  \right\},
\]
which is isomorphic to
\begin{align}\label{fanw}
\left\{x\in \mathbb{R}^{2n}_{\ge0}: \sum_{k=2a_i-1}^{2a_{i+1}-2}x_k=a_{i+1}-a_{i},  \sum_{k\in[2j]+2a_i-2}x_k\ge j ,\ \text{for}\ 1\le j\le a_{i+1}-a_i-1  \right\}.
\end{align}
Notice that in \eqref{fanw}, if $ a_{i+1}-a_i\le j\le n-1$, since $x_{2a_{i+1}-1}+\cdots+x_{2a_{i+1}-2+2k}\ge k$ for any $1\le k\le n-1$, we see that the inequality
\[
\sum_{k\in[2j]+2a_i-2}x_k=a_{i+1}-a_i+\sum_{k=2a_{i+1}-1}^{2j+2a_i-2}x_k\ge a_{i+1}-a_i+ (j-a_{i+1}+s_i)=j
\]
automatically holds.
That is to say,  $\mathcal{P}(C_{a_{i+1}-a_i})$ is isomorphic to
\begin{align*}
\left\{x\in \mathbb{R}^{2n}_{\ge0}: \sum_{k=2a_i-1}^{2a_{i+1}-2}x_k=a_{i+1}-a_{i},  \sum_{k\in[2j]+2a_i-2}x_k\ge j ,\ \text{for}\ 1\le j\le n-1  \right\}.
\end{align*}

Similarly, 
\begin{align*}
  \mathcal{P}(C_{n-a_{r}+a_1})=\left\{x\in \mathbb{R}^{2n}_{\ge0}: \sum_{k\in[2(n-a_{r}+a_1)]}x_k=n-a_{r}+a_1,  \sum_{k\in[2j]}x_k\ge j ,\ \text{for}\ 1\le j\le n-a_{r}+a_1-1  \right\}.
\end{align*}
After shifting all the indices by $2a_r-2$, we see that
$\mathcal{P}(C_{n-a_{r}+a_1})$ is isomorphic to 
\begin{align}\label{115sw}
\left\{x\in \mathbb{R}^{2n}_{\ge0}: \sum_{2a_{r}-1}^{2n+2a_1-2}x_k=n-a_{r}+a_{1},  \sum_{k\in[2j]+2a_r-2}x_k\ge j ,\ \text{for}\ 1\le j\le n-a_{r}+a_{1}-1 \right\}.
\end{align}
Notice that in \eqref{115sw}, if $ n-a_r+a_1\le j\le n-1$, since $x_{2a_{1}-1}+\cdots+x_{2a_{1}-2+2k}\ge k$ for any $1\le k\le n-1$, we see that the inequality
\[
\sum_{k\in[2j]+2a_r-2}x_k=n-a_r+a_1+\sum_{k=2a_{1}-1}^{2j-2n+2a_r-2}x_k\ge n-a_r+a_1+ (j-n+a_r-a_1)=j
\]
always holds. Thus $\mathcal{P}(C_{n-a_{r}+a_1})$ is isomorphic to
\begin{align*} 
\left\{x\in \mathbb{R}^{2n}_{\ge0}: \sum_{k=1}^{2a_1-2}x_k+\sum_{2a_{r}-1}^{2n}x_k=n-a_{r}+a_{1},  \sum_{k\in[2j]+2a_r-2}x_k\ge j ,\ \text{for}\ 1\le j\le n-1 \right\}.
\end{align*}

Therefore, the set of nonnegative solutions of equation \eqref{spl} is  isomorphic to the Cartesian product of the Catalan matroid polytopes $\mathcal{P}(C_{n-a_{r}+a_1})$ and $\mathcal{P}(C_{a_{i+1}-a_i})$ for $1\le i\le r-1$. That is,
 \begin{align}\label{Q}
Q(A)= \mathcal{P}(C_{a_{2}-a_1})\times \mathcal{P}(C_{a_{3}-a_2})\times\cdots\times \mathcal{P}(C_{a_{r}-a_{r-1}})\times\mathcal{P}(C_{n-a_{r}+a_1}).
\end{align}
Therefore,
 \begin{align*}
\indicator{Q(A)}= 
\indicator{C_{a_{2}-a_1}}\cdots \indicator{C_{a_{r}-a_{r-1}}}\cdot \indicator{C_{n-a_{r}+a_1}}=\indicator{C_\lambda}, 
\end{align*}
as required.\qed

For a partition  $\lambda=(\lambda_1,\ldots,\lambda_k)$   of 
$n$, let $\ell(\lambda)$ denote the number of parts of $\lambda$, i.e., $\ell(\lambda)=k$. We can also write $\lambda=(1^{\alpha_1},\ldots,n^{\alpha_n})$, where $\alpha_i$ is the number of appearances of $i$ in $\lambda$. Denote $p(\lambda)$ by the number of permutations of the parts of $\lambda$, that is,
\[p(\lambda)=\frac{\ell(\lambda)!}{\alpha_1!\alpha_2!\cdots\alpha_n!}.\]
For a composition $\gamma=(\gamma_1,\ldots,\gamma_k)$ of $n$,
if there exists an integer $m>0$  such that  $\gamma_i=\gamma_{i+m}$  for all $1\le i\le k-m$,  then we call $m$ a period of $\gamma$. Let $T(\gamma)$ denote the least period of $\gamma$. If there does not exist such a positive integer $m$, set $T(\gamma)=k$.  For instance, if $\gamma=(2,2,1,2,2,1,2,2,1)$, then $T(\gamma)=3.$ If $\gamma=(2,1,3,2)$, then $T(\gamma)=4.$

\begin{lem}\label{lm2}
Let $\lambda\vdash n$ be a partition. Then we have
\[
\#\{A\subseteq[n]\colon \lambda(A)=\lambda\}=\frac{n}{\ell(\lambda)}p(\lambda).
\]
\end{lem}

\pf Let $\lambda=(\lambda_1,\ldots,\lambda_k)$.  Given a permutation $\gamma=(\lambda_{i_1},\lambda_{i_2},\ldots,\lambda_{i_k})$ of $\lambda$, we can construct $n$  subsets $A\subseteq[n]$ such that $\lambda(A)=\lambda$ as below
\begin{align}\label{S1}
\{j,j+\lambda_{i_1},j+\lambda_{i_1}+\lambda_{i_2},\ldots,j+\lambda_{i_1}+\cdots+\lambda_{i_{k-1}}\},
\end{align}
where  $1\le j\le n$ and  any  $\lambda_m$ is understood as $\lambda_{m-n}$ if $m>n$. It is easy to see that all subsets   $A\subseteq[n]$ such that $\lambda(A)=\lambda$ can be generated in this way. 
However, there are repetitions in \eqref{S1}. To be exact, it is easy to see that  the set
\[
\{j',j'+\lambda_{i_1},j'+\lambda_{i_1}+\lambda_{i_2},\ldots,j'+\lambda_{i_2}+\cdots+\lambda_{i_{k-1}}\},
\]
where $j'=j+\lambda_{i_1}+\cdots+\lambda_{i_{T(\gamma)}}$ is the same as the set in \eqref{S1}. Therefore, for a given a permutation $\gamma$ of $\lambda$,  each subset $A$ of $[n]$  such that $\lambda(A)=\lambda$ is generated $\frac{\ell(\lambda)}{T(\gamma)}$  times.

On the other hand,  different permutations of $\lambda$ may generate the same subset $A$. More precisely, for a shifting  $\gamma'=(\lambda_{i_2},\lambda_{i_3},\ldots,\lambda_{i_k},\lambda_{i_1})$ of $\gamma$, if $\gamma'\neq\gamma$, then we can construct the following subset
\begin{align}\label{S2}
\{j'',j''+\lambda_{i_2},j''+\lambda_{i_2}+\lambda_{i_3},\ldots,j''+\lambda_{i_2}+\cdots+\lambda_{i_{k}}\},
\end{align}
where $j''=j+\lambda_{i_1}$. Since $j''+\lambda_{i_2}+\cdots+\lambda_{i_{k}}=j+n$, the two sets in \eqref{S1} and \eqref{S2} are the same. Each permutation $\gamma$ of $\lambda$ can shift $T(\gamma)$ different times, all generating the same subset. Thus the total number of repetitions of each $A$ is $\ell(\lambda)$. 

Since there are $p(\lambda)$ permutations of $\lambda$, and for each permutation  of $\lambda$, we  constructed $n$ subsets $A$ as in \eqref{S1}, but there are $\ell(\lambda)$ copies of $A$ constructed in this way, we conclude that the total number of subsets $A$ with $\lambda(A)=\lambda$  is $\frac{n}{\ell(\lambda)}p(\lambda)$. 
\qed

Suppose that
$\tau=(\tau_1,\ldots,\tau_k)$  and $\sigma=(\sigma_1,\ldots,\sigma_s)$ are  two compositions of $n$. we say that $\tau$ and $\sigma$ are  {\it equivalent }, denoted as $\tau\sim \sigma$, if they have the same multiset of parts. Obviously, if $\tau\sim\sigma$,  then $\indicator{U_{\tau},t}=\indicator{U_{\sigma},t}$. We say that $\tau$ is a {\it refinement } of $\sigma$, denoted as $\tau\prec\sigma$, if 
 each $\sigma_i$ can be obtained by adding a sequence of consecutive factors  of $\tau$.    More precisely,  there exists $i_1,i_2,\ldots,i_{s-1}$  such that 
\[\sigma_1=\tau_1+\cdots+\tau_{i_1},
\sigma_2=\tau_{i_1+1}+\cdots+\tau_{i_2},\ldots,\sigma_s=\tau_{i_{s-1}+1}+\cdots+\tau_{k}.\] Let $F(\lambda)$ be the set of refinements of $\lambda$. For instance, let $\lambda=(3,2)$, then $F(\lambda)=\{(3,2),(2,1,2),(1,1,1,2),$ $(3,1,1),(2,1,1,1),(1,1,1,1,1)\}$.  

For a partition $\lambda\vdash n$, let $K_\lambda$ denote the set of permutations on $[n]$ with {\it cycle type} $\lambda$. Recall that a permutation $w\in S_n$ has cycle type $\lambda=(1^{\alpha_1},\ldots,n^{\alpha_n})$ means that when converted into cycles, $w$ has $\alpha_i$ cycles of length $i$, for $i=1,2,\ldots,n$. It is well known that
\[|K_{\lambda}|=\frac{n!}{z_{\lambda}},\]
 see, for example, Sagan \cite{Sg}.

\begin{prop}\label{Theoindicator}
 In the space of indicator functions,
\[
\indicator{\mathcal{P}(C_n)} = \sum_{\lambda \vdash n} \frac{1}{z_\lambda} \, \indicator{\mathcal{P}(U_\lambda)}.
\]
Moreover, the Catalan point  is an interior point in the polytope of all matroids $\overline{\Omega}_{r,n}$.
\end{prop}
\pf
 By \eqref{incexc}, Lemma \ref{lm1} and Lemma \ref{lm2}, we have
\begin{align}		
\indicator{U_n}=\sum_{\lambda\vdash n}(-1)^{\ell(\lambda)-1}\frac{np(\lambda)}{\ell(\lambda)}\indicator{C_{\lambda}}.
\end{align}
Then
 \begin{align}\label{expcn}
\indicator{C_n}=\frac{1}{n}\indicator{U_n}+\sum_{\lambda\vdash n,\lambda\ne(n)}(-1)^{\ell(\lambda)}\frac{p(\lambda)}{\ell(\lambda)}\indicator{C_{\lambda}}.
\end{align}
Since  $z_{(n)}=\frac{1}{n}$ and $\indicator{U_{(n)}}=\indicator{U_n}$, by  \eqref{expcn}, it suffices to show that
\begin{align}\label{gds}
\sum_{\lambda\vdash n,\lambda\ne(n)}(-1)^{\ell(\lambda)}\frac{p(\lambda)}{\ell(\lambda)}\indicator{C_{\lambda}}=
\sum_{\lambda\vdash n,\lambda\ne(n)} \frac{1}{z_{\lambda}}\indicator{U_{\lambda}}.
\end{align}

We make induction on $n$.
By induction,  
\begin{align}\label{compp}
&\sum_{\lambda\vdash n,\lambda\ne(n)}(-1)^{\ell(\lambda)}\frac{p(\lambda)}{\ell(\lambda)}\indicator{C_{\lambda}}\nonumber\\[5pt]
&=\sum_{\lambda\vdash n,\lambda\ne(n)}(-1)^{\ell(\lambda)}\frac{p(\lambda)}{\ell(\lambda)}\prod_{j=1}^{\ell(\lambda)}\left(\sum_{\sigma^j\vdash \lambda_j} \frac{1}{z_{\sigma^j}} \indicator{U_{\sigma^j}}\right)\nonumber\\[5pt]
&=\sum_{\lambda\vdash n,\lambda\ne (n)}(-1)^{\ell(\lambda)}\frac{p(\lambda)}{\ell(\lambda)}\sum_{\substack{\sigma\in F(\lambda)}} \frac{1}{z_{\sigma^1}z_{\sigma^2}\cdots z_{\sigma^{\ell(\lambda)}}} \indicator{U_{\sigma}},
\end{align}
where, given a refinement $\sigma\in F(\lambda)$ of $\lambda$, there is   a unique way to cut  $\sigma$ into $\ell(\lambda)$ factors  $\sigma^1,\sigma^2,\ldots,\sigma^{\ell(\lambda)}$ such that each $\sigma^i\vdash\lambda_i$ is a partition of $\lambda_i$.

In order to prove \eqref{gds}, we need to show that for any  partition $\mu$ of $n$ with $\mu\neq(n)$, the coefficient of $f(U_{\mu},t)$ in \eqref{compp}  is equal to $1/z_\mu$. Since  $\indicator{U_{\sigma},t}=\indicator{U_{\mu},t}$ if $\sigma\sim\mu$,   the  coefficient of $\indicator{U_{\mu},t}$ in \eqref{compp} is
\begin{align*}  
 \sum_{\substack{\lambda\vdash n,\lambda\neq(n)}}(-1)^{\ell(\lambda)}\frac{p(\lambda)}{\ell(\lambda)} \sum_{\sigma\in F(\lambda),\sigma\sim\mu } \frac{1}{z_{\sigma^1}z_{\sigma^2}\cdots z_{\sigma^{\ell(\lambda)}}}.
\end{align*}
That is, we need to show that
\begin{align*}  
 \sum_{\substack{\lambda\vdash n,\lambda\neq(n)}}(-1)^{\ell(\lambda)}\frac{p(\lambda)}{\ell(\lambda)} \sum_{\sigma\in F(\lambda),\sigma\sim\mu } \frac{1}{z_{\sigma^1}z_{\sigma^2}\cdots z_{\sigma^{\ell(\lambda)}}}=\frac{1}{z_\mu},
\end{align*}
 which is equivalent to
\begin{align*} 
 \sum_{\substack{\lambda\vdash n,\lambda\neq(n)}}(-1)^{\ell(\lambda)}
 \frac{p(\lambda)}{\ell(\lambda)} 
 \sum_{\sigma\in F(\lambda),\sigma\sim\mu } \frac{\lambda_1!}{z_{\sigma^1}}\frac{\lambda_2!}{z_{\sigma^2}}\cdots \frac{\lambda_{{\ell(\lambda)}}!}{z_{\sigma^{{\ell(\lambda)}}}}\cdot \binom{n}{\lambda_1,\ldots,\lambda_{\ell(\lambda)}}=\frac{n!}{z_{\mu}}.
\end{align*}
Therefore, it is enough to show that for any  partition $\mu=(\mu_1,\ldots,\mu_r)\vdash n, \mu\neq(n)$, there holds
\begin{align} \label{laststep}
 \sum_{\substack{\lambda\vdash n,\lambda\neq(n)}}(-1)^{\ell(\lambda)}
 \frac{p(\lambda)}{\ell(\lambda)} 
 \sum_{\sigma\in F(\lambda),\sigma\sim\mu } |K_{\sigma^1}||K_{\sigma^2}|\cdots |K_{\sigma^{\ell(\lambda)}}|\cdot \binom{n}{\lambda_1,\ldots,\lambda_{\ell(\lambda)}}=|K_{\mu}|.
\end{align}

Suppose that    $w\in K_{\mu}$, that is, when converted in cycles, $w$  has cycle type $\mu$. In particular, $w$ has $r$ cycles. We need to show that the total number of appearances of $w$ in the left-hand side of \eqref{laststep} is equal to 1.
To this end, we discuss according to the number of parts of $\lambda$ in the left-hand side of \eqref{laststep}. 
We claim that for any fixed $k=\ell(\lambda)\ge 2$,   $w$ is counted in 
\begin{align}\label{cyctype}
\sum_{\lambda\vdash n,\ell(\lambda)=k}(-1)^k\frac{p(\lambda)}{k}
 \sum_{\sigma\in F(\lambda),\sigma\sim\mu } |K_{\sigma^1}||K_{\sigma^2}|\cdots |K_{\sigma^{k}}|\cdot \binom{n}{\lambda_1,\ldots,\lambda_{k}}
\end{align}
by $(-1)^k(k-1)!S(r,k)$ times, where $S(r,k)$ is Stirling number of the second kind.

Let $\lambda=(\lambda_1,\ldots,\lambda_k)=(1^{\alpha_1},\ldots,n^{\alpha_n})\vdash n$. Notice that if $\sigma\in F(\lambda)$, then  $\sigma$ admits a unique decomposition into factors $\sigma^1,\ldots,\sigma^k$ such that $\sigma^{i}\vdash \lambda_i$. Choose an ordered  partition of $[n]$, i.e.,
\begin{align}\label{sets}
T_1=\{j_1,\ldots,j_{\lambda_1}\},
T_2=\{j_{\lambda_1+1},\ldots,j_{\lambda_1+\lambda_2}\},
\ldots,
T_k=\{j_{s},\ldots,j_{n}\},
\end{align}
where $s=n-(\lambda_1+\cdots+\lambda_{k-1})+1$, such that $|T_i|=\lambda_i$, $T_i\cap T_j=\emptyset$ for $i\neq j$ and $T_1\cup\cdots\cup T_k=[n]$. There are $\binom{n}{\lambda_1,\ldots,\lambda_{k}}$ such ordered partitions of $[n]$.  
For $1\le i\le k$,  let  $T_i$ form a permutation of cycle type $\sigma^i$, there are $|K_{\sigma^i}|$ ways.
Suppose that $w$ is counted in \eqref{cyctype} at least once. Then each $T_i$ in \eqref{sets} must be a union of some cycles of $w$. Moreover, although each $T_i$ can form $|K_{\sigma^i}|$ different permutations of cycle type $\sigma^i$, there is exactly one way to obtain the cycles of $w$. Thus
$w$ is counted exactly once by a given ordered partition of $[n]$ in \eqref{sets}.

Assume that $|T_i|=|T_j|$ for some $i<j$ and $T_i$ forms a permutation of cycle type $\sigma^i$, $T_j$ forms a permutation of cycle type $\sigma^j$. Then we can 
exchange the positions of $T_i$  and $T_j$  in \eqref{sets} to obtain    $k$ ordered   disjoint subsets of $[n]$: 
\[T_1,\ldots, T_j,\ldots,T_i,\ldots,T_k,\]
and let them form $k$ permutations of cycle types $\sigma^1,\ldots,\sigma^j,\ldots,\sigma^i,\ldots,\sigma^k$, respectively. Clearly, $w$ can also be obtained once in this way. In other words, for a given $\lambda\vdash n$ of length $k$,   $w$  is counted  in
\[
\sum_{\sigma\in F(\lambda),\sigma\sim\mu } |K_{\sigma^1}||K_{\sigma^2}|\cdots |K_{\sigma^{k}}|\cdot \binom{n}{\lambda_1,\ldots,\lambda_{k}}
\]
by
$\alpha_1!\alpha_2!\cdots\alpha_n!=k!/p(\lambda)$ times. 

To prove the claim, we  need to sum over all partitions $\lambda\vdash n$ with $\ell(\lambda)=k$. If $w$  is counted once by \eqref{cyctype}, then each  $T_i$ in \eqref{sets} is a union of some cycles of $w$. In other words, for a given $w\in K_{\mu}$, we can merge the $r$ cycles of $w$ to obtain $k$ sets $T_i$  in \eqref{sets}. If we  run over all possible ways of ``cycles merging'' of $w$, then we actually   obtain all partitions of $n$ with exactly $k$ parts. Obviously,   there are $S(r,k)$ such ways. Therefore,    $w$ is counted  in   \eqref{cyctype} by
\[
(-1)^k\frac{p(\lambda)}{k}\frac{k!}{p(\lambda)}S(r,k)=
(-1)^k(k-1)!S(r,k)
\]
times.

By summing for all $k\ge2$ and noticing that
$S(r,k)=kS(r-1,k)+S(r-1,k-1),$
we conclude that the total number of $w\in K_{\mu}$ counted in the left-hand side of \eqref{laststep} is
\begin{align*}
&\sum_{k=2}^{r}(-1)^k(k-1)!S(r,k)\\
&\quad=\sum_{k=2}^{r}(-1)^k(k!S(r-1,k)+(k-1)!S(r,k-1))\\
&\quad=S(r-1,1)+(-1)^nn!S(r-1,n)+\sum_{k=2}^{r-1}(k!S(r-1,k)-k!S(r-1,k))\\
&\quad=S(r-1,1)\\
&\quad=1,
\end{align*}
as desired.
 \qed

For example, let $\mu=(2,2,1,1)$ and $w=(1,4)(2,5)(3)(6)\in K_{\mu}$. We enumerate how many times $w$ is counted in \eqref{cyctype}. The set of partitions $\lambda\vdash 6, \lambda\neq(6)$ such that 
there exists a refinement $\sigma\in F(\lambda)$ with $\mu\sim\sigma$ is
\[
\{(2,2,1,1),(3,2,1),(2,2,2),(4,1,1),(5,1),(4,2),(3,3)\}.
\]

For $\lambda=(2,2,1,1)$,  among the $\binom{6}{2,2,1,1}=180$ possible ordered partitions of $\{1,\ldots,6\}$, there are four that can be used to obtain $w$, i.e.,  
\begin{align*}
T_1&=\{1,4\}, T_2=\{2,5\}, T_3=\{6\}, T_4=\{3\},\\
T_1&=\{1,4\}, T_2=\{2,5\}, T_3=\{3\}, T_4=\{6\},\\
T_1&=\{2,5\}, T_2=\{1,4\}, T_3=\{6\}, T_4=\{3\},\\
T_1&=\{2,5\}, T_2=\{1,4\}, T_3=\{3\}, T_4=\{6\}.
\end{align*}
For each of ordered partition above, we let $T_1,T_2$ form permutations with cycle type $\sigma^1=\sigma^2=(2)$ and let $T_3,T_4$ form permutations with cycle type $\sigma^3=\sigma^4=(1)$.
Thus  $w$ is counted $(-1)^{\ell(\lambda)}\frac {p(\lambda)} {\ell(\lambda)}\cdot4 = 6$ times in \eqref{cyctype} for $\lambda = (2,2,1,1)$.

For $\lambda=(3,2,1)$, there are four ordered partitions of $\{1,\ldots,6\}$  that can be used to obtain $w$, i.e.,
\begin{align*}
T_1&=\{1,4,3\}, T_2=\{2,5\}, T_3=\{6\},\\
T_1&=\{1,4,6\}, T_2=\{2,5\}, T_3=\{3\},\\
T_1&=\{2,5,3\}, T_2=\{1,4\}, T_3=\{6\},\\
T_1&=\{2,5,6\}, T_2=\{1,4\}, T_3=\{3\}.
\end{align*}
For each of them, let $T_1$ form a permutation of cycle type $\sigma^1=(2,1)$,  $T_2,T_3$ form   permutations of cycle types $\sigma^2=(2), \sigma^3=(1)$, respectively. Thus $w$ is counted $(-1)^{\ell(\lambda)}\frac{p(\lambda)}{\ell(\lambda)}\cdot4=-8$ times in \eqref{cyctype} for $\lambda=(3,2,1)$. 

Similarly,   $w$ is  counted $-2$ times for $\lambda=(4,1,1)$, $-2$ times for $\lambda=(2,2,2)$, $3$ times for $\lambda=(4,2)$, $2$ times for $\lambda=(5,1)$, and $2$ times for $\lambda=(3,3)$, respectively. In total, $w$ is counted $6-8-2-2+3+2+2=1$ time.

\begin{prop}\label{coroindicatorab}
\[
\indicator{\mathcal{P}(C^{a,b}_n)} = \sum_{\lambda \vdash n} \frac{1}{z_\lambda} \, \indicator{\mathcal{P}(U^{a,b}_\lambda)}.
\]
Moreover, the $(a,b)$-Catalan point  is an interior point in the polytope of all matroids $\overline{\Omega}_{r,n}$.
\end{prop}
Since the proofs of Proposition \ref{Theoindicator} and Proposition \ref{coroindicatorab} are nearly identical, we omit the proof of Proposition \ref{coroindicatorab} for notational simplicity.


In fact, using Ferroni and Fink's terminology of the "matroid polytope" – 
where uniform matroids are extremal (see \cite[Example 4.4]{F-allmatroids}), 
and more generally direct sums of uniform matroids are also vertices, see Bonin~\cite[Theorem~8.1]{Bonin2025}. 
Theorem~\ref{Theo_nonextrme} provides a reformulation of our main result.

\begin{re}
\label{re:other-Schubert}
The non-extremality phenomenon of $(a,b)$-Catalan matroids depends on a strong cyclic symmetry. To be specific, the matroid polytope $\mathcal{P}(C_n^{a,b})$ maps to an isomorphic copy of itself under a cyclic coordinate transformation. Moreover, the collection of such copies produced by all cyclic shifts precisely forms a subdivision of a hyperplane. This structural property is crucial for our argument.

By contrast, for general Schubert matroids, while cyclic (or other) coordinate transformations may still yield isomorphic copies of the base polytope, the union of these copies typically does not assemble into a hyperplane subdivision. For example, for a Schubert matroid $\mathrm{SM}_n(S)$ with $r(S) = (a, b, c, d, a, b, c, d)$, where $a,b,c,d$ are positive integers, cyclic transformations indeed generate isomorphic polytopes. But the union of any such collection does not fill a hyperplane in the ambient space. This lack of tiling prevents the extension of our approach.
Therefore, the non-extremality phenomena for $(a,b)$-Catalan matroids do not generally extend to arbitrary Schubert matroids.
\end{re}

\section{Valuative Invariants of Catalan Matroids}
Proposition~\ref{coroindicatorab} immediately yields a closed formula for arbitrary valuative invariants $f$ on $(a,b)$-Catalan matroids, expressed in terms of direct sums of uniform matroids, as formalized in Corollary~\ref{CoroVI}.

\subsection{The Ehrhart Polynomials}

 In \cite{FL}, Fan and Li provided a formula for the Ehrhart polynomials of Schubert matroids. As an application, they obtained a recursive formula for the Ehrhart polynomials  of $(a,b)$-Catalan matroids, and conjectured that  $(a,b)$-Catalan matroids are Ehrhart positive.  

Although Ferroni \cite{Fer3} showed that not all matroids are Ehrhart positive, disproving  conjectures of  De Loera, Haws and K\"{o}ppe \cite{Lorea} and  Castillo and Liu \cite{Liu},  various specific subclasses of matroids turn out to be Ehrhart positive. For instance,
uniform matroids, minimal matroids, sparse paving Schubert matroids, panhandle matroids, and all matroids of rank two  are  Ehrhart positive, see \cite{FL,Fer,Fer2,Fer4,panhandleposi}.  In this note, we express the Ehrhart polynomials of $(a,b)$-Catalan matroids  as   positive combinations of   Ehrhart polynomials of uniform matroids, and thus obtain the Ehrhart positivity of  $(a,b)$-Catalan matroids.

Given a polytope $\P$ and a positive integer $t$, the $t$-dilation of $\P$ is defined as $t\P=\{t\alpha\colon\alpha\in \P\}$. Let $i(\P,t)=|t\P\cap \mathbb{Z}^n|$ denote the number of lattice points in $t\P$. For integral polytopes, $i(\P, t)$ is a polynomial in $t$, called the Ehrhart polynomial of $\P$. For simplicity, we write $i(M, t)$ for $i(\P(M), t)$.

\begin{coro} \label{coroEhr} Let $a,b$ be positive integers, we have
\begin{align}
i(C_n^{a,b},t)=\sum_{\lambda\vdash n} \frac{1}{z_{\lambda}} i(U_{\lambda}^{a,b},t).
\end{align}
\end{coro}
Since Ferroni\cite{Fer} has shown that uniform matroids are Ehrhart positive, we obtain the Ehrhart positivity of   $(a,b)$-Catalan matroids.

Recall that notched rectangle matroids are Schubert matroids $\SM_n(S)$ with $r(S)=(a,b,c,d)$. In particular,  $r(S)=(a,b,c,1)$ corresponds to the  panhandle matroid ${\texttt{Pan}}_{b+1,a+b,a+b+c+1}$, which was conjectured to be Ehrhart positive in \cite{panhandle}, and proved in \cite{panhandleposi}. When  $c=1$,  $r(S)=(a,b,1,1)$ corresponds to the panhandle matroid  ${\texttt{Pan}}_{b+1,a+b,a+b+2}$.

 For simplicity,   write $i(r(S), t)$ for $i(\SM_n(S), t)$.
By \cite[Theorem 1.5]{FL},  there holds 
\begin{align*}
i((a,b,a,b),t)&=\frac{1}{2}i(U_{2b,2a+2b},t)
+\frac{1}{2}i(U_{b,a+b},t)^2,
\\
i((a,a,b,b),t)&=\frac{1}{2}i(U_{a+b,2a+2b},t)
+\frac{1}{2}i(U_{a,2a},t)i(U_{b,2b},t),\\
i((1,1,a,a+1),t)&=\frac{1}{2}(t+2)i(U_{a+1,2a+2},t).
\end{align*}
In this paper, we express $i(C_n^{a,b},t)$ in terms of the Ehrhart polynomials of uniform matroids. 
Therefore, it seems  feasible to express other families of Ehrhart polynomials in terms of the Ehrhart polynomials of uniform matroids. McGinnis \cite{Mc} showed that ${\texttt{Pan}}_{b+1,a+b,a+b+2}$ is Ehrhart positive by using intricate combinatorial arguments.   We  propose  the 
following Conjecture, which manifestly shows  the Ehrhart positivity of ${\texttt{Pan}}_{b+1,a+b,a+b+2}$.

\begin{conj}
    Let $a,b$ be positive integers, we have
\[
i({\texttt{Pan}}_{b+1,a+b,a+b+2},t)=i((1,1,a,b),t)=\left(\frac{a+1}{a+b+1}t+1\right)\cdot i(U_{b,a+b+1},t).\]
\end{conj}

We remark that the Schubert matroid  $\SM_n(S)$ is  paving  if and only if  $r(S)=(a,1,b,c)$. In particular, $\SM_n(S)$ is a paving panhandle matroid if and only if  $r(S)=(1,1,a,b)$.

As an extension of the preceding discussion, we are now at a position to consider the volume of $\mathcal{P}(C_n^{a,b})$.
The following corollary follows easily from Corollary \ref{coroEhr}.

\begin{coro}\label{volume}The volume of $\mathcal{P}(C^{a,b}_n)$ is given by    \begin{align}\label{volabC}        \operatorname{vol}(\mathcal{P}(C_n^{a,b}))=\frac{1}{n(2n-1)!}A(n(a+b)-1,nb),        \end{align}where $A(n,k)$ denotes the Eulerian number, which counts the number of permutations on $[n]$ with exactly $k-1$ descents.\end{coro}
\subsection{The Tutte Polynomials}


Let $M$ be a matroid on $E$. The Tutte polynomial of a matroid $M$ is defined as the bivariate polynomial $T_{M}(x, y) \in \mathbb{Z}[x, y]$ given by
$$
T_{M}(x, y)=\sum_{A \subseteq E}(x-1)^{\mathrm{rk}(E)-\mathrm{rk}(A)}(y-1)^{|A|-\mathrm{rk}(A)}.
$$

\begin{exam}
For $1 \leqslant k \leqslant n-1$, the Tutte polynomial of the uniform matroid $U_{k, n}$ is given by
$$
T_{U_{k, n}}(x, y)=\sum_{i=1}^k\binom{n-i-1}{n-k-1} x^i+\sum_{i=1}^{n-k}\binom{n-i-1}{k-1} y^i.
$$
\end{exam} 


The Tutte polynomial is a valuative invariant for matroids, as established in earlier works of    Speyer \cite{Speyer} and    Ardila,  Fink, and Rinc\'{o}n \cite{AFR}. Ferroni and Schr\"{o}ter\cite{F-VI} provides an alternative concise proof of this property using convolution techniques.

As a corollary of Corollary~\ref{CoroVI}, we provide a closed formula of Tutte polynomials for arbitrary $(a,b)$-Catalan matroids.

\begin{coro} \label{coroTut} Let $a,b$ be positive integers, we have
\begin{align}
T_{C_n^{a,b}}(x,y)=\sum_{\lambda\vdash n} \frac{1}{z_{\lambda}} T_{U_\lambda^{a,b}}(x,y).
\end{align}
\end{coro}

\subsection{The Kazhdan-Lusztig polynomials}

\subsubsection{The KL polynomials}

Elias, Proudfoot, and Wakefield \cite{kl-of-matroid} introduced the Kazhdan–Lusztig polynomial $P_M$  for loopless matroids $M$, establishing an analogy with the Kazhdan–Lusztig theory for Coxeter groups. This polynomial is defined recursively using the lattice of flats. Specifically, they provided a recurrence relation for the coefficients of its polynomial expansion:
\[
P_{U_{k,n}}(t) = \sum_i c^i_{k,n}\,t^i.
\]
Gao, Xie, Yang, and Zhang \cite{Yang} later  derived an explicit closed‑form expression for these coefficients.

\begin{theo}[\cite{Yang}, Theorem 1.3]\label{thm-ukl3}
    For any $n$, $k$, and $0\leq i\leq \lfloor\frac{k-1}{2}\rfloor$, we have
    \begin{align}\label{ukl3}
    c_{k,n}^i={\frac{1}{k-i} \binom{n}{i}
    {\sum _{h=0}^{n-k-1}\binom{k-i+h}{h+i+1}\binom{i-1+h}{h} }}.
    \end{align}
\end{theo}

As a direct consequence of Corollary~\ref{CoroVI}, we obtain:
\begin{coro}  Let $a,b$ be positive integers, we have
\begin{align}\label{corokl}
P_{C_n^{a,b}}(t)=\sum_{\lambda\vdash n} \frac{1}{z_{\lambda}} P_{U_\lambda^{a,b}}(t).
\end{align}
\end{coro}

\begin{exam}

\[
P_{\mathrm{C}_n}(t)=P_{\mathrm{C}_n^{1,1}}(t)=
\begin{cases}
1, & n=2,\\[4pt]
3t+1, & n=3,\\[4pt]
15t+1, & n=4,\\[4pt]
45t^2+55t+1, & n=5,\\[4pt]
473t^2+185t+1, & n=6,\\[4pt]
1092t^3+3239t^2+612t+1, & n=7.\\[4pt]
\end{cases}
\]

\end{exam}

\subsubsection{The inverse KL polynomial}
Gao and Xie \cite{inverse-kl} derived an explicit  expression for the inverse Kazhdan–Lusztig polynomial of the uniform matroid $U_{k,n}$.
\begin{theo}[\cite{inverse-kl},Theorem 3.3] \label{inverse-kl-uniform}
	For any uniform matroid $U_{k,n}$ with $n,k\geq 1$, we have
	$$Q_{U_{k,n}}(t)=\binom{n}{k}\sum_{j=0}^{\lfloor (k-1) /2\rfloor}\frac{(n-k)(k-2j)}{(n-k+j)(n-j)}\binom{k}{j}t^j.$$
\end{theo}
 
\begin{coro} \label{coroinv-kl} Let $a,b$ be positive integers, we have
\begin{align}
Q_{C_n^{a,b}}(x,y)=\sum_{\lambda\vdash n} \frac{1}{z_{\lambda}} Q_{U_\lambda^{a,b}}(x,y).
\end{align}
\end{coro}

\begin{exam}
    \[
Q_{\mathrm{C}_n}(t)=Q_{\mathrm{C}_n^{1,1}}(t)
\begin{cases}
2, & n=2,\\[4pt]
3t+5, & n=3,\\[4pt]
19t+14, & n=4,\\[4pt]
45t^2+92t+42, & n=5,\\[4pt]
396t^2+405t+132, & n=6,\\[4pt]
1092t^3+2491t^2+1705t+429, & n=7.\\[4pt]
\end{cases}
\]

\end{exam}

\subsubsection{$Z$-polynomial}
Proudfoot, Xu, and Young \cite{proudfoot2017z} defined the $Z$-polynomial for any matroid  as follows
\[
Z_M(t)\;=\;\sum_{F\in \mathcal{L}(M)}t^{\rk(M_F)}\,P_{M^F}(t),
\]
where $\mathcal{L}(M)$ denotes the lattice of flats, $M^F$ represents the contraction of $M$ by the flat $F$, and $P_{M^F}(t)$ the Kazhdan–Lusztig polynomial associated with $M^F$.  In particular, for the uniform matroid $U_{k,n}$ they proved the identity
\begin{equation}\label{kltoz}
Z_{U_{k,n}}(t)
= t^k \;+\;\sum_{i=1}^{k}\binom{n}{n-k+i}\,t^{\,k-i}\,P_{U_{i,n-k+i}}(t).
\end{equation}
From \eqref{kltoz} it follows that $Z_{U_{k,k+1}}(t)$ coincides  with the Narayana polynomial.  Expressing $Z_{U_{k,n}}(t)$ as
\[
Z_{U_{k,n}}(t)\;=\;\sum_{i\ge0}z_{k,n}^i\,t^i,
\]
Gao et al. \cite{Yang} provided an explicit formula for the coefficients $z_{k,n}^i$.

\begin{theo}[\cite{Yang}, Theorem 1.5]\label{thm-uz1}
For any $k$, $n$, and $0\leq i\leq k$,
we have
\begin{align}
z_{k,n}^i=\frac{\binom{n}{i+n-k} \binom{n}{i}}{\binom{n}{n-k}}
{\sum _{h=0}^{n-k-1}  \frac{i (h-n+k+1)+n-k}{(h+1) (n-k)} \binom{i-1+h}{h} \binom{k-i+h}{h}}\label{uz1}.
\end{align}
\end{theo}

We derive the corresponding formula for the Catalan case:
\begin{coro}
    Let $a,b$ be positive integers, for the $(a,b)$-Catalan matroid $C_n^{a,b}$, we have
\begin{align}
Z_{C_n^{a,b}}(t)=\sum_{\lambda\vdash n} \frac{1}{z_{\lambda}} Z_{U_\lambda^{a,b}}(t).
\end{align}
\end{coro}

\begin{exam}
    \[
Z_{\mathrm{C}_n}(t)=Z_{\mathrm{C}_n^{1,1}}(t)=
\begin{cases}
t^2+3t+1, & n=2,\\[4pt]
t^3+8t^2+8t+1, & n=3,\\[4pt]
t^4+22t^3+50t^2+22t+1, & n=4,\\[4pt]
t^5+64t^4+278t^3+278t^2+64t+1, & n=5,\\[4pt]
t^6+196t^5+1433t^4+2619t^3+1433t^2+196t+1, & n=6,\\[4pt]
t^7+625t^6+7010t^5+20596t^4+\cdots+1, & n=7.\\[4pt]
\end{cases}
\]

\end{exam}

\subsection{Whitney polynomial}
An important invariant of a matroid $M$ of rank $k$ is given by its Whitney numbers of the second kind, namely
\[
W_i \;=\;\#\{\text{flats of $M$ of rank $i$}\}, i=0,1,\dots,k.
\]
Following Mason \cite{Mason}, the rank generating function of the lattice of flats $\mathcal L(M)$ is the Whitney polynomial of a rank $k$ matroid $M$
\[
W_M(x)\;=\;\sum_{i=0}^k W_i\,x^i.
\]
By the valuative properties of the flag $f$-vector, $W_M(x)$ defines a valuative invariant of $M$, and we already know the following formula.

\begin{exam}
    The Whitney polynomial of $U_{k,n}$ is
    \[
W_{U_{k,n}}(x) = \sum_{i=0}^{k-1} \binom{n}{i}\,x^i + x^k.
\]
\end{exam}

Then we can derive a closed  expression for the Whitney polynomial of any $(a,b)$‑Catalan matroid.

\begin{coro}\label{coroWhitney}
Let $a,b$ be positive integers, the Whitney polynomial of the $(a,b)$-Catalan matroid $C_n^{a,b}$  is as follows
\[
W_{C_n^{a,b}}(x)
=\sum_{\lambda\vdash n} \frac{1}{z_{\lambda}}W_{U_\lambda^{a,b}}(x).
\]
\end{coro}

\begin{exam}
    \[ 
    W_{\mathrm{C}_n}(t)=
W_{\mathrm{C}_n^{1,1}}(t)=
\begin{cases}
t^2+3t+1, & n=2,\\[4pt]
t^3+8t^2+5t+1, & n=3,\\[4pt]
t^4+22t^3+19t^2+7t+1, & n=4,\\[4pt]
t^5+64t^4+67t^3+34t^2+9t+1, & n=5,\\[4pt]
t^6+196t^5+232t^4+144t^3+53t^2+11t+1, & n=6,\\[4pt]
t^7+625t^6+804t^5+573t^4+261t^3+76t^2+13t+1, & n=7.\\[4pt]
\end{cases}
\]

\end{exam}

\section*{Acknowledgement}The authors are grateful to Neil J.Y. Fan and Luis Ferroni for helpful conversations. YC and MY are partially supported by NSF grant No. 12471314.

\small

Y.M. C{\scriptsize HEN}, D{\scriptsize EPARTMENT OF} M{\scriptsize ATHEMATICS}, S{\scriptsize ICHUAN} U{\scriptsize NIVERSITY}, C{\scriptsize HENGDU} 610064, P.R. C{\scriptsize HINA.}

\vspace{-.2cm}
{Email address:  \tt \href{mailto:ym\_chen@stu.scu.edu.cn}{ym\_chen@stu.scu.edu.cn}}

Y. L{\scriptsize I}, D{\scriptsize EPARTMENT OF} M{\scriptsize ATHEMATICS}, S{\scriptsize ICHUAN} U{\scriptsize NIVERSITY}, C{\scriptsize HENGDU} 610064, P.R. C{\scriptsize HINA.}

\vspace{-.2cm}{ Email address: \tt\href{mailto:liyaao@stu.pku.edu.cn}{liyaao@stu.pku.edu.cn}}

M. Y{\scriptsize AO}, D{\scriptsize EPARTMENT OF} M{\scriptsize ATHEMATICS}, S{\scriptsize ICHUAN} U{\scriptsize NIVERSITY}, C{\scriptsize HENGDU} 610064, P.R. C{\scriptsize HINA.} 

\vspace{-.2cm}{ Email address:  \tt \href{mailto:yaom@stu.scu.edu.cn}{yaom@stu.scu.edu.cn}}

\end{document}